\documentclass[12pt]{article}
\usepackage{amsfonts}
\usepackage{amsmath}
\usepackage{amssymb}
\usepackage{latexsym}
\usepackage{amscd}
\usepackage{bbm}
\usepackage{url}
\usepackage[left=3.2cm,right=3.2cm, top=2.5cm,bottom=2.5cm,bindingoffset=0cm]{geometry}

\newenvironment{Proof}{\par\noindent{\it Proof.}} {\hfill$\scriptstyle\blacksquare$}

\newcommand\address{\noindent\leavevmode
	\medskip
	\noindent
	Kateryna Tatarko, \\
	Dept.~of Math.~and Stat.~Sciences,\\
	University of Alberta, \\
	Edmonton, Alberta, Canada, T6G 2G1.\\
	\texttt{\small
		e-mail:  tatarko@ualberta.ca}
}
\date{}

\newcounter{theorem}[section]
\newtheorem{theorem}{Theorem}
\newtheorem{lemma}[theorem]{Lemma}
\newtheorem{cor}[theorem]{Corollary}
\newtheorem{defi}[theorem]{Definition}
\newtheorem{prop}[theorem]{Proposition}

\newtheorem{Remark}[theorem]{Remark}

\newcommand{\R}{\mathbb{R}}
\newcommand{\norm}[1]{\|{#1}\|}

\newcommand{\RomanNumeralCaps}[1]
{\MakeUppercase{\romannumeral #1}}

\DeclareRobustCommand{\rchi}{{\mathpalette\irchi\relax}}
\newcommand{\irchi}[2]{\raisebox{\depth}{$#1\chi$}} 

\title{An upper bound on the smallest singular value of a square random matrix}
\author{
	Kateryna Tatarko
}

\date{}

\begin{document}
 \hangindent=1.5cm
\hangafter=-1 \noindent

\maketitle

\begin{abstract}
	Let $A = (a_{ij})$ be a square $n\times n$ matrix with i.i.d. zero mean and unit variance entries. It was shown by Rudelson and Vershynin in 2008 that 
	the upper bound for the smallest singular value $s_n(A)$ is of order $n^{-\frac12}$ with probability close to one under the additional assumption that the entries of $A$ satisfy $\mathbb{E}a^4_{11} < \infty$.
	 We remove  the  assumption  on the fourth moment and show the upper bound assuming only $\mathbb{E}a^2_{11} = 1.$
\end{abstract}



\section{Introduction}

The extremal singular values have been attracting the attention of scientists in different disciplines such as mathematical physics  or geometric functional analysis. In particular, they play an important role in numerical analysis as the condition number, which is the ratio of the largest to the smallest singular value, is a measure for the worst-case loss of precision in a computational problem. Much is known about the behavior of the largest singular value and we refer the reader to~\cite{BSY, YBK}.
The study of the behavior of the smallest singular value goes back to von Neumann and his collaborators concerning numerical inversion of large matrices, where they conjectured (see~\cite{vonN,vonN1}) that the smallest singular value is of order $n^{-\frac12}$ with probability close to one. Estimates of similar type for the case of Gaussian matrices (i.e., matrices with i.i.d. standard normal entries) were obtained by Edelman in~\cite{E} and Szarek in~\cite{S}. For estimates on extremal singular values which were acquired while studying the problem of the approximation of covariance matrices, we refer to \cite{ALPT,GLPT, MP, T2}. Various bounds for the smallest singular value have been obtained under rather weak assumptions on the rows of the matrix in \cite{KM, O, Y, Y1}.  For lower bounds on the smallest singular value of random matrices with independent but not identically distributed entries see a recent result by Cook~\cite{C}.

Rudelson and Vershynin in~\cite{RV1,RV2,RV3} studied the behavior of the smallest singular value of matrices with i.i.d. subgaussian entries. They showed (see~\cite{RV1,RV2}) that the smallest singular value of a square random matrix $A$ with i.i.d. subgaussian entries is of order $n^{-\frac12}.$  In particular, in~\cite{RV1} they proved that for given $t\geq 2 $ there are $C >0\text{ and } u\in(0, 1)$ depending only on the subgaussian moment of entries of $A$ such that $$\mathbb{P} \left( s_n(A)~>~t n^{-\frac12}\right)~\leq~C\frac{\log{t}}{t} + u^n.$$
 Nguyen and Vu in~\cite{NV} showed an exponential bound for the above probability, which improves the linear bound by Rudelson and Vershynin. A lower bound for rectangular subgaussian matrices was obtained  in~\cite{RV3}. A recent result of Wei (see~\cite{W}) provides  upper bounds on intermediate singular values of rectangular matrices with subgaussian entries. The corresponding lower bounds were obtained in~\cite{R1}.
 
Recently, in~\cite{RT} a new technique was developed, which allowed Rebrova and Tikhomirov to prove a lower bound for $s_n(A)$ of square matrices of order $n^{-\frac12}$ under the assumption that the L\'evy concentration function of entries of $A$ is bounded. Namely, they showed the small ball probability estimate: $$\forall \varepsilon> 0: \quad \mathbb{P}\left( s_n(A) \leq \varepsilon n^{-\frac12} \right)~\leq~C\varepsilon + u^n,$$ where $C > 0 \text{ and } u \in (0, 1)$ depend only on the law of $a_{11}.$ Notice that any random variable $\xi$ with $\mathbb{E}\xi~=~0 \text{ and } \mathbb{E}\xi^2 = 1$ has a bounded L\'evy concentration function, therefore the above statement is valid for matrices with assumptions only on the second moment of entries.

The goal of this note is to show that the upper bound on the smallest singular value holds for square matrices with heavy-tailed entries. We prove the following theorem.
\begin{theorem}\label{main theorem}
	Let $A = (a_{ij})$ be an $n\times n$ matrix whose entries are i.i.d. random variables with $\mathbb{E}a_{ij} = 0$ and $\mathbb{E}a_{ij}^2 = 1$.
	Then there exists an absolute constant $C > 0$ 
	such that for every~$\varepsilon > 0$
	\begin{equation*}
	\mathbb{P}\left(s_n(A) > \frac1{\varepsilon^2} n^{-\frac12} \right) \leq C\varepsilon + \frac{C}{\sqrt{n}}. 
	\end{equation*} 
\end{theorem}
We expect that the dependence on $\varepsilon$ can be improved to ${\varepsilon}^{-1},$ but our proof gives only~${\varepsilon^{-2}}.$

We now briefly describe the ideas of proof of Theorem~\ref{main theorem}. 

To estimate the smallest singular value of a random matrix $A$ we will use the following equivalence, which holds for every $\lambda \geq 0,$
$$
s_n(A) \leq \lambda \quad \Longleftrightarrow \quad \exists x\in S^{n - 1}: \ \norm{Ax} \leq \lambda.
$$
We will show that there exists $x \in \R^n$ such that $\norm{x} \leq \tau$ and $\norm{A^{-1}x} \geq \eta \sqrt{n}$ for some $\tau, \eta > 0,$ which implies $s_n(A)\leq \frac{\tau}{\eta\sqrt{n}}.$ Let us describe the main difficulty in our proof. It is well-known that $A^{-1}x$ behaves differently depending on the structure of $x.$ We follow~\cite{LPRT, LPRT2} and roughly speaking split the unit sphere into two parts consisting of vectors of small dimensions and vectors with bounded $\ell_\infty$ norm. To deal with vectors of the second type, we use ideas introduced in~\cite{RV2}, namely we use the essential least common denominator (see the definition below). Denote by $B$ the transpose of the first $n - 2$ columns of matrix $A$. To show that the essential least common denominator of vectors in the null space of a matrix $B$ has exponential decay with high probability, in~\cite{RV1} the authors used a standard $\varepsilon$-net argument, namely, for a given $\varepsilon$-net $\mathcal{N}$ on a subset $S\subset S^{n - 1}$ one has $$\inf\limits_{y \in S}\norm{By}~\geq~\inf\limits_{y' \in \mathcal{N}}\Big(\norm{By'}~-~\norm{B}\norm{y - y'}\Big).$$ This procedure relies on an upper bound for the operator norm $\norm{B}$, which is of order $n^{\frac12}$  with exponentially high probability under the subgaussian moment assumption on the entries of $B.$ Moreover, as can be seen in~\cite{G, YBK}, 
one has that $\norm{B} \leq C\sqrt{n}$ under the assumption of bounded fourth moments (see \cite{Lat,LvHY} for independent but not identically distributed entries). However, in the settings of Theorem~\ref{main theorem}, it is not guaranteed that the operator norm $\norm{B}$ has a good upper bound. Moreover, if the fourth moment is unbounded, it is known that $\frac{\norm{B}}{\sqrt{n}} \to \infty$ \big(\cite{BSY, Sil, YBK}, see also \cite{LS} for quantitative estimates\big). To overcome this difficulty, we use a recent technique developed by Rebrova and Tikhomirov in~\cite{RT}. Starting with a standard $\varepsilon$-net on $S\subset S^{n - 1}$ 
 we construct a new net on $S$ which is a $\left(C\varepsilon \sqrt{n}\right)$-net with respect to the pseudometric $\norm{B(x - y)}$ with probability close to one. This allows us to circumvent the use of the operator norm $\norm{B}.$ 


\section{Preliminaries}

By $e_1, \dots, e_n$ we denote the canonical basis of $\R^n$ equipped with the canonical inner product $\left<\cdot,\cdot\right>$ and Euclidean norm $\norm{\cdot}.$

Let $A$ be an $n \times n$ matrix with real entries. Then the singular values $s_j(A),$ $ j\leq n$, of the matrix $A$ are the eigenvalues of $\sqrt{A^{*}A}$, which are arranged in non-increasing order: $s_1(A) \ge s_2(A) \ge \dots \ge s_n(A).$  In particular, the largest and the smallest singular values  $s_1(A)$ and $s_n(A)$ are
\begin{equation*}
s_1(A) =  \|A\| \quad \text{ and } \quad s_n(A) = \frac{1}{\|A^{-1}\|},
\end{equation*}

\noindent where $\| A \| = \sup\limits_{\|x\| = 1}\|Ax\| $ is the operator norm from $ \ell_2^n$ to $\ell_2^n$, 
and $A^{-1}$ is the inverse from the image of $A$. 

Recall that for a given metric space $X$, an $\varepsilon$-net $\mathcal{N}$ in $X$ is a subset of $X$ such that any point of $X$ is within distance at most $\varepsilon$ from points of $\mathcal{N}.$ 

A system $\left(E_k,F_k\right)_{k = 1}^n$ of vectors $(E_k)_{k = 1}^n$ and $(F_k)_{k = 1}^n$ in an $n$-dimensional Hilbert space $H$ is called a biorthogonal system if $\left<E_k,F_s\right> = \delta_{k,s}$ for all $k,s \in \{1, \dots, n\}$. The system is called complete if it spans the entire space $H.$ The next proposition contains some well-known properties of biorthogonal systems (see~\cite{RV1}, Proposition~2.1).

\newpage
\begin{prop}\leavevmode 
	\begin{itemize}
		\item[\textup{(i)}] Let $\left(E_k\right)_{k = 1}^n$ be a linearly independent system of vectors in an $n$-dimensional Hilbert space $H.$ Then there exist unique vectors $\left(F_k\right)_{k = 1}^n$ such that $\left(E_k,F_k\right)_{k = 1}^n$ is a complete biorthogonal system in $H$.
		\item[\textup{(ii)}] If $\left(E_k,F_k\right)_{k = 1}^n$ is a complete biorthogonal system in $H$, then $$\norm{F_k} = \frac{1}{\textup{dist}(E_k, H_k)} \text{ \ \ \ for } k~=~1, \dots, n,$$
		where $H_k = \textup{span}(E_i)_{i\ne k}.$
		\item[\textup{(iii)}] If $A$ is an $n\times n$ invertible matrix, then $\left(Ae_k, \left(A^{-1}\right)^te_k\right)_{k = 1}^n$ is a complete biorthogonal system. 
	\end{itemize} \label{biorthogonal systems}	
\end{prop}

We will also need the notion of the so-called L\'evy concentration function of a random variable $\xi$,
\begin{equation*}
\mathcal{L}(\xi, \varepsilon) = \sup_{\lambda \in \mathbb{R}} \mathbb{P}\left( \left|\xi - \lambda\right| \leq \varepsilon\right), \ \varepsilon \geq 0.
\end{equation*}
\noindent In other words, the L\'evy concentration function measures how likely a random variable $\xi$ enters a small ball in the space.
As we mentioned above, any random variable $\xi$ with $\mathbb{E}\xi = 0$ and $\mathbb{E}\xi^2 = 1$ satisfies the condition $$\mathcal{L}(\xi, v)  \le u$$ for some constants $u\in(0, 1)$ and $v >0$ determined by the law of $\xi$. Therefore, we don't add this constraint to the formulation of our main result Theorem~\ref{main theorem}, but state it only in terms of finiteness of the second moment of entries.

In order to find an upper bound for the smallest singular value $s_n(A)$, we will consider a partition of the sphere into sets of compressible and incompressible vectors. Such an idea to split the sphere into two parts and to use an estimate involving the norm of a matrix, in order to bound the smallest singular value first appeared in~\cite{LPRT} and was formalized later (see~\cite{RV2}) in the following definition.

\begin{defi}
	Let $\delta, \rho \in (0, 1).$ A vector $x \in \mathbb{R}^n$ is called $\left(\delta n\right)${\it-sparse} if $$|\textup{supp}(x)| < \delta n.$$ A vector $x \in S^{n - 1}$ is called {\it compressible} if $x$ is within Euclidean distance $\rho$ from the set of all $\delta n$-sparse vectors. Otherwise, a vector $x \in S^{n - 1}$ is called {\it incompressible}. 
\end{defi}
The sets of compressible and incompressible vectors will be denoted by $${\rm Comp} = {\rm Comp}_n(\delta, \rho) \quad \quad \text{and} \quad \quad {\rm Incomp} = {\rm Incomp}_n(\delta, \rho),$$ respectively.

Since the set of compressible vectors is essentially of the smaller dimension, the following simple result shows that one can find an $\varepsilon$-net on the set of compressible vectors $\text{Comp}$ with small cardinality. 

\begin{lemma}\label{euclid_net}
For any $\delta,\rho \in (0,\,1]$ a set of  compressible unit vectors $\textup{Comp}_n(\delta, \rho)$ admits a $\left(2\rho\right)$-net $\mathcal{N}$ of $\textup{Comp}_n(\delta, \rho)$ of  cardinality $$|\mathcal{N}| \leq \left(\frac{e}{\delta}\right)^{\delta n}\Big(\frac{5}{\rho}\Big)^{\delta n}.$$
\end{lemma}
\begin{Proof}
By definition, for every $x\in \text{Comp}_n(\delta, \rho)$ there exist $x'\in S^{n - 1}$ such that $|\text{supp}(x')|\leq\delta n $  and $ ||x - x'|| \leq \rho.$ Thus, to find a $\left(2\rho\right)$-net on a set of compressible vectors, it is enough to find a Euclidean $\rho$-net on the set of sparse vectors. For a fixed coordinate subspace of dimension $\delta n$, the standard volumetric estimate gives a $\rho$-net of a cardinality at most $(1 + \frac2{\rho})^{\delta n}.$ Applying a union bound over all coordinate subspaces, we have that the set of compressible vectors $\text{Comp}_n(\delta, \rho)$ admits an Euclidean $\left(2\rho\right)$-net of cardinality
\begin{equation*}
|\mathcal{N}|\leq \dbinom{n}{\delta n}\left( 1+\frac{2}{\rho}\right)^{\delta n} \leq 
\left( \frac{e}{\delta}\right) ^{\delta n} \left( \frac{5}{\rho}\right) ^{\delta n}.
\end{equation*}
\end{Proof}

We will need a couple of results from~\cite{RT}. The following theorem allows us to refine a given $\varepsilon$-net $\mathcal{N}$ on a subset of the unit sphere to an $\left(\frac{\varepsilon C}{\delta}\sqrt{n}\right)$-net $\widetilde{\mathcal{N}}$ on the same subset of the sphere with respect to pseudometric $\norm{A(x  -y)} $ with high probability.
\begin{theorem}[\cite{RT}, Theorem $A^\star$]
	Let $\delta \in (0, \frac14]$, $\varepsilon \in (0,\frac12],$ $n\geq \frac{1}{4\delta}$, $S\subset S^{n - 1}$ be a subset of the sphere, and $\mathcal{N}\subset S$ be an $\varepsilon$-net on $S$ in the Euclidean metric. Then there exists a deterministic subset $\widetilde{\mathcal{N}}\subset S$ with $$|\widetilde{\mathcal{N}}|\leq \textup{exp}\left(13\delta n \ln{\frac{2e}{\delta}}\right)|\mathcal{N}|$$ such that for an $n\times n$ random matrix $A$ with i.i.d. zero mean and unit variance entries, with probability at least $1- 4 \textup{exp}(-\frac{\delta n}{8}),$ the set $\widetilde{\mathcal{N}}$ is an $\left(\frac{\varepsilon C }{\delta}\sqrt{n}\right)$-net on $S$ with respect to the pseudometric $\textup{d}(x,y)=\|A(x-y)\|$, where $x,y \in S^{n-1}.$ \label{refine net}
\end{theorem}

\begin{Remark}
	One can check that Theorem~\ref{refine net} holds for a $(n - 2)\times n$ matrix~$A$. 
\end{Remark}

The next lemma gives a strong probability estimate for a fixed unit vector.

\begin{lemma}[\cite{RT}, Lemma 4.9] \label{lemma2}
	Let $\xi$ be a random variable with $\mathcal{L}(\xi,\tilde{v}) \leq \tilde{u}$ for some $\tilde{v}>0 \text{ and } \tilde{u}\in(0, 1).$ Then there are $v>0$ and $u\in(0, 1)$ depending only on $\tilde{u}, \tilde{v}$ such that for an $(n - 2)\times n$   random matrix $A$ with i.i.d. entries equidistributed with $\xi$ and for any $y\in S^{n - 1}$ one has
	\begin{equation*}
	\mathbb{P}\left(\|Ay\|\leq v\sqrt{n}\right)\leq u^{n - 2}.
	\end{equation*}	
\end{lemma}

In order to obtain the small ball probability estimate for a random sum, we need the notion of the essential least common denominator. It measures the closeness of the scaled vector $x \in \R^n$  to $\mathbb{Z}^n$. This notion was introduced in~\cite{RV2,RV3} (see also~\cite{TV}) and for more detailed description see~\cite{R}.

\begin{defi}
	For parameters $\alpha > 0$  and $r \in (0, 1)$, {\it the essential least common denominator of a vector} $x \in \R^n$ is defined as
	\begin{equation*}
	\textup{LCD}_{\alpha, r} (x) = \inf\left\{t > 0: \textup{dist}(tx, \mathbb{Z}^n) < \min\left( r\norm{tx}, \alpha\right) \right\}.
	\end{equation*}
Then {\it the essential least common denominator of a subspace} $H \subset \R^n$ is defined as
	\begin{equation*}
	\textup{LCD}_{\alpha, r} (H) = \inf\left\{ \textup{LCD}_{\alpha, r}(x): x\in H, \norm{x} = 1 \right\}.
	\end{equation*}
\end{defi}
Later we will use this definition with a small constant $r$, and  a small multiple $\alpha$ of $\sqrt{n}.$

The next result gives a small ball probability estimate of a random sum. It is essentially Theorem 3.4 in~\cite{RV1}.

\begin{theorem}\label{small ball prob}
	Let $u\in (0,1).$ Let $\xi_1, \dots, \xi_n$ be i.i.d. zero mean random variables such that $\mathcal{L}(\xi_1, 1) \leq u$  and $x = (x_1, \dots, x_n) \in S^{n - 1}.$ Then, for every $\alpha > 0,$ $r \in (0, 1)$ and for every $\varepsilon > 0$
	one has 
	$$
	\mathcal{L}\left( \sum_{i=1}^{n} x_i\xi_i, \varepsilon\right)  \leq \frac{C}{r\sqrt{1 - u}}\left(\varepsilon + \frac1{\textup{LCD}_{\alpha, r}(x)} \right) + Ce^{-2\alpha^2(1 - u)},
	$$
	where $C>0$ is an absolute constant.
\end{theorem}

In words, the theorem provides useful upper bounds on the small ball probability which depend on the additive structure of the coefficients $x_1, \dots, x_n.$ The less structure the coefficients carry, the more spread the distribution of a random sum is, and the less the small ball probability is.


\section{Proof of the Theorem~\ref{main theorem}}

To prove the boundedness of the smallest singular value of the type
\begin{equation*}
s_n(A) \leq Ln^{-\frac12},
\end{equation*}
where $L>0$, it is enough to show that there exists $x\in\R^n$ such that $\norm{x} \leq \tau$ and $\norm{A^{-1}x} \geq \eta n^{-\frac12}$ for some $\tau,\eta >0.$ 

We follow the ideas of Rudelson and Vershynin in~\cite{RV1}. Consider the columns $X_i = Ae_i$ of a matrix $A$ and the rows $\widetilde{X}_i = \left(A^{-1}\right)^t e_i$ of an inverse matrix $A^{-1}.$ Let $H_i$ denote the span of all column vectors except the $i$-th, i.e.
$$H_i = \text{span}\left(X_1, \dots,X_{i - 1},X_{i+1},\dots, X_n\right),$$
and $H_{i,j}$ denote the span of all column vectors except the $i$-th and $j$-th ($i < j$), i.e.
$$H_{i, j} = \text{span}\left(X_1, \dots,X_{i - 1},X_{i+1},\dots, X_{j - 1}, X_{j + 1}, \dots, X_n\right).$$
Let $P_1$ denote the orthogonal projection in $\R^n$ onto the subspace $H_1$ and let
$$
x = X_1 - P_1X_1.
$$
Then $x$ is orthogonal to $H_1$. Since our matrix $A$ is invertible and  $\text{dim}\,\text{ker}P_1 = 1,$ then we also have that $\norm{x} = \text{dist}\left( X_1, H_1\right).$

Note that by Markov's inequality, we have
\begin{equation}
\mathbb{P}\left(\norm{x} > \tau \right) \leq \frac{\mathbb{E} \norm{x}^2}{\tau^2}, \quad \quad \tau>0. \label{expectation of x}
\end{equation}

Let $f_n$ be a normal vector of the $(n-1)$-dimensional subspace $H_1.$ Then, the vector $x$ can be represented as $x =  \left<X_1, f_n \right>f_n,$ and the norm of $x$ is $$\norm{x} = \left| \left<X_1, f_n \right> \right| = \left|\sum_{i = 1}^{n}a_{i1}f_n^i  \right|.$$ Hence,
\begin{equation}\label{isotropic}
\mathbb{E}\left|\sum_{i = 1}^{n}a_{i1}f_n^i  \right|^2 = \mathbb{E}\left(\sum_{i = 1}^{n}a_{i1}^2\left(f_n^i\right)^2 +  \sum_{i \ne j }a_{i1}a_{j1}f_n^if_n^j \right) = \sum_{i = 1}^{n}\left(f_n^i\right)^2 \mathbb{E}a_{i1}^2 = 1
\end{equation}
(this fact also follows from the fact that vector $X_1$ is isotropic). 
Then by~\eqref{expectation of x},
$$
\mathbb{P}\left(\norm{x} > \tau \right) \leq \frac{1}{\tau^2}, \quad \quad \tau >0. 
$$

Now we estimate $\norm{A^{-1}x}.$ Note that 
$$
\norm{A^{-1}x} = \norm{A^{-1}X_1 - A^{-1}P_1X_1} = \norm{e_1 - A^{-1}P_1Ae_1}.
$$
Since the vector $P_1Ae_1$ belongs to  $\text{span}\left\{Ae_2, \dots,Ae_n\right\}$, then $A^{-1}P_1Ae_1$ is orthogonal to $e_1.$ Therefore, using $P_1\widetilde{X}_1 = 0$ and denoting $Y_k = P_1\widetilde{X}_k,$ $\ k~\in~\{2,\dots,n\}$, we obtain
\begin{align}\label{norm_inequality}
\norm{A^{-1}x}^2 &= \norm{e_1}^2 + \norm{A^{-1}P_1X_1}^2 > \norm{A^{-1}P_1X_1}^2 = \sum_{k = 1}^{n} \left<A^{-1}P_1X_1, e_k \right>^2 \nonumber \\  &= 
\sum_{k = 1}^{n} \left<X_1, P_1\widetilde{X}_k \right>^2 
= \sum_{k = 2}^{n} \left<X_1, Y_k \right>^2.
\end{align}

\bigskip

The following lemma provides the relation between families of vectors $(X_k)_{k = 2}^n$ and $(Y_k)_{k = 2}^n$. 
\begin{lemma}[\cite{RV1}, Lemma 2.1]
If $\left(X_k, Y_k\right)_{k = 2}^n$ is defined as above, then it is a complete biorthogonal system in $H_1.$\label{biorthogonal_system}
\end{lemma}

The following is a  consequence of the uniqueness in Proposition~\ref{biorthogonal systems} (i).
\begin{cor} \label{corollary}
	The system of vectors $(Y_k)_{k = 2}^n$ is uniquely determined by the system $(X_k)_{k = 2}^n$. In particular, the system $(Y_k)_{k = 2}^n$ and the vector $X_1$ are independent.
 \end{cor}

By Proposition~\ref{biorthogonal systems} (ii), we have $\norm{Y_k} = \frac{1}{\text{dist}\left(X_k, H_{1,k}\right)}.$ Therefore, we can rewrite~\eqref{norm_inequality} as
\begin{equation}
\norm{A^{-1}x}^2 \geq \sum_{k = 2}^{n}\frac{1}{1/\norm{Y_k}^2}\big<\frac{Y_k}{\norm{Y_k}}, X_1\big>^2 = \sum_{k = 2}^{n}\left(\frac{a_k}{b_k}\right)^2, \label{sum}
\end{equation}
where 
\begin{equation}\label{coefficients}
a_k =\left| \big<\frac{Y_k}{\norm{Y_k}}, X_1\big> \right| \  \text{ and } \  b_k = \frac{1}{\norm{Y_k}} = \text{dist}\left(X_k, H_{1, k}\right).
\end{equation}
This reduces our problem to bounding $a_k$ from above and $b_k$ from below. Without loss of generality, we can do it for $k = 2,$ since the same argument carries over to any $k\in\left\{2, \dots, n\right\}$

We split the unit sphere into sets of compressible and incompressible vectors. Our next goal is to show that the orthogonal complement $H^\perp_{1, 2}$ consists of incompressible vectors with high probability. Consider an~$(n - 2)~\times~n$~matrix $B$ with columns $X_3, \dots, X_n.$ Since the subspace $H_{1, 2}$ is the span of the independent random vectors $X_3, \dots, X_n,$ we have $H^\perp_{1, 2} \subset \text{ker}(B).$ We want to show:
\begin{equation}
\forall x\in \text{Comp}: \quad \norm{Bx} \geq \lambda \sqrt{n}, \label{compressible}
\end{equation}
that is, with high probability compressible vectors do not belong to the kernel of matrix $B$ (the parameter $\lambda$ will be determined later).

To deal with compressible vectors, we need the following proposition, which is essentially Proposition 5.2 from~\cite{RT}, where it was proved for $n\times n$ matrices. For the sake of completeness, we provide the proof for $(n - 2) \times n$ matrices.

\begin{prop} \label{prob_estimate_copm}
	Let $\xi$ be a centered random variable with unit variance such that $\mathcal{L}(\xi, \tilde{v}) \leq \tilde{u}$ for some $\tilde{v} > 0$ and $\tilde{u} \in (0,\,1).$ Let $n \in \mathbb{N}$ and let $\Gamma$ be an $(n - 2)\times n$ random matrix with i.i.d. entries equidistributed with $\xi.$ Then there are numbers $\theta,\, v >0$ and $u \in (0,\,1)$ depending only on $\tilde{u}, \tilde{v}$ such that for $\textup{Comp} = \textup{Comp}_{n}(\rho, \rho)$ we have
	$$
	\mathbb{P}\left( \inf_{x \in \textup{Comp}}\norm{\Gamma y} < v\sqrt{n}\right) \leq 5u^{n - 2}.
	$$ 
\end{prop}
\begin{Proof}	
The main idea of the proof is to apply the union bound over the set of compressible vectors $\text{Comp}.$ 
In Theorem~\ref{refine net} take $\delta \in (0, \frac14]$ such that $$e^{13n\delta \ln{\frac{2e}{\delta}}} \leq u^{-\frac {n - 2}3}.$$ Then define the parameter $\rho \in (0, \frac16]$ in such a way that
$$
\left(\frac{5e}{\rho^2}\right)^{\rho n} \leq u^{-\frac{n - 2}{3}} \ \ \text{ and } \ \  \frac{3\rho C}{\delta} \leq \frac v2,
$$
where $C > 0$ is a universal constant taken from Theorem~\ref{refine net}.

By Lemma~\ref{euclid_net}, there is a Euclidean $\left(2\rho\right)$-net $\mathcal{N} \subset \text{Comp}$ of cardinality $$|\mathcal{N}| \leq \left(\frac{5e}{\rho^2}\right)^{\rho n}~\leq~u^{-\frac{n - 2}{3}}.$$ Now we refine this net using Theorem~\ref{refine net}, and as a result with probability at least $1-4\text{exp}\left(-\frac{\delta n}{8}\right)$ we obtain a $\left(\frac{2\rho C}{\delta}\sqrt{n}\right)$-net $\widetilde{\mathcal{N}} \subset \text{Comp}$ with respect to the pseudometric $\norm{\Gamma(x - y)}$ which has cardinality $|\mathcal{\widetilde{N}}| \leq u^{-2(n - 2)/3}.$ In other words, for every $x \in \text{Comp}$ there exists $x' = x'(x) \in \text{Comp}$ such that $$\norm{\Gamma(x - x')} \leq \frac{2\rho C}{\delta}\sqrt{n} \leq \frac v2 \sqrt{n}.$$

Applying the union bound over $\mathcal{\widetilde{N}}$ to the relation from Lemma~\ref{lemma2}, we get
\begin{equation*}
\mathbb{P}\left(\norm{\Gamma x'} < v \sqrt{n} \text{ for some } x' \in \widetilde{\mathcal{N}}\right) \leq |\widetilde{\mathcal{N}}|u^{n - 2} \leq u^{\frac {n  -2}3}.
\end{equation*}
On the other hand, the construction of $\widetilde{\mathcal{N}}$ implies that
\begin{equation*}
\mathbb{P}\left( \inf_{x \in \text{Comp}} \norm{\Gamma x} < \inf_{x' \in \widetilde{\mathcal{N}}} \norm{\Gamma x'} - \frac v2 \sqrt{n} \right) \leq 4\text{exp}\left(-\frac{\delta n}{8}\right).
\end{equation*}
Therefore,
\begin{equation*}
\mathbb{P}\left(\norm{\Gamma x} < \frac v2 \sqrt{n} \text{ for some } x \in \text{Comp}\right) \leq u^{\frac {n - 2}3} + 4\text{exp}\left(-\frac{\delta (n - 2)}{8}\right).
\end{equation*}
Taking the maximum of $u^{\frac13}$ and $e^{-\frac{\delta}8}$ gives the desired result.
\end{Proof}


\bigskip

The next proposition states that the least common denominator of any incompressible vector in $\mathbb{R}^n$ is of order at least $\sqrt{n}.$ This proposition is Lemma~$6.1$ from~\cite{R} (note that the proof does not depend on the parameter $\alpha$).
\begin{prop}
	For any parameters $\theta, \rho \in (0, 1)$ there are parameters $r, \gamma > 0$ such that for every $\alpha > 0$ any vector $x \in \textup{Incomp}_n(\theta, \rho)$ satisfies $\textup{LCD}_{\alpha, r}(x) \geq \gamma \sqrt{n}.$ \label{incomp_of_order_root_n}
\end{prop}

Recall that $B$ is a $(n - 2)\times n$ matrix with columns $X_3, \dots, X_n.$ Since $H^\perp_{1, 2} \subset~\text{ker}(B)$, for the set $\text{Comp}$ the following implication holds:
$$
\text{if } \quad \inf\limits_{y \in \text{Comp}} \norm{By} > 0 \quad \text{ then } \quad H^\perp_{1, 2}\cap \text{Comp} = \emptyset. 
$$
Applying Proposition~\ref{prob_estimate_copm} to the matrix $B$, we get 
$$
\mathbb{P}\left(\inf\limits_{y \in \text{Comp}} \norm{By} \geq v\sqrt{n}\right) \geq 1 - 5u^{n - 2}.
$$
Therefore,  $H^\perp_{1, 2} \cap\text{Comp} = \emptyset$ with probability at least $1 - 5u^{n - 2}$, or in other words,
$$
\mathbb{P}\left(H^\perp_{1, 2} \cap S^{n - 1} \subseteq \text{Incomp}\right) \geq 1 - 5u^{n - 2},
$$
which means that the subspace $H^\perp_{1, 2}$ consist of incompressible vectors with probability close to one.
By Proposition~\ref{incomp_of_order_root_n}, we obtain that for some $u \in (0, 1)$
\begin{equation}
\mathbb{P}\left( \text{LCD}_{\alpha, r}\left(H^\perp_{1, 2}\right) \geq \gamma \sqrt{n} \right) \geq 1 - 5u^{n - 2}, \label{incompressible}
\end{equation}
where $\alpha$ is a small multiple of $\sqrt{n}.$

Recall that the coefficients $a_k$ and $b_k$ were introduced in~\eqref{coefficients}. To ensure that the lower bound for $b_2$ is satisfied with high probability, we condition on $H_{1, 2}$ and use Markov's inequality and the fact that $X_2$ is isotropic (see \eqref{isotropic}). More precisely, we obtain
\begin{equation}
\mathbb{P}\left(b_2 = \text{dist}\left(X_2, H_{1, 2} \right)\geq t \,|\, H_{1,2} \right) \leq \frac{\mathbb{E}\,\text{dist}\left(X_2, H_{1, 2} \right)^2}{t^2} \leq \frac{2}{t^2}, \quad t>0. \label{b2 bound}
\end{equation}

Let $\mathcal{E} = \left\{ \text{LCD}_{\alpha, r}\left( H_{1, 2 }^\perp \right) \geq \gamma \sqrt{n} \text{ and } b_2 < t \,|\, H_{1,2} \right\}.$ Combining the two estimates~\eqref{incompressible} and~\eqref{b2 bound}, we get that 
\begin{equation}
\mathbb{P}\left( \mathcal{E} \right) \geq 1 - \frac2{t^2} - 5u^{n - 2}. \label{event}
\end{equation}

Since we conditioned on the subspace $H_{1, 2}$, we may fix a realization of vectors $\left(X_j\right)_{j = 2}^n$ for which the statement~\eqref{event} holds. Thus by the uniqueness in Corollary~\ref{corollary} the vector $Y_2$ is also fixed. For convenience, we further consider the normalized vector $Y = \frac{Y_2}{\norm{Y_2}}.$ By Lemma~\ref{biorthogonal_system} we know that $\left(X_j\right)_{j = 2}^n$ and $\left(Y_j\right)_{j = 2}^n$ form a biorthogonal system, in particular $Y$ is orthogonal to $\left(X_j\right)_{j = 3}^n.$ Thus, $Y \in H_{1, 2}^\perp$. Since the event $\mathcal{E}$ in~\eqref{event} holds, we know that 
$$
\text{LCD}_{\alpha, r} (Y) \geq \gamma \sqrt{n}.
$$ 
Now we proceed to bound the coefficient $a_2.$ Recall that
$$
a_2 = |\left<Y, X_1\right>| =| \sum_{k = 1}^{n} Y^kX_1^k |$$
and $Y^k$ are coefficients such that $\sum_{k = 1}^{n} \left( Y^k\right)^2 = 1$ and $X_1^k$ are i.i.d. random variables with zero mean and $\mathcal{L}\left(X_1^k, 1\right) < u$ for some $u\in(0, 1).$ Applying Theorem~\ref{small ball prob} with $\alpha = c\sqrt{n}$ for some small absolute constant $c>0,$  we obtain for $\varepsilon > 0,$ $u \in (0, 1)$ and $r \in (0, 1)$
\begin{eqnarray}
\mathbb{P}_{X_1}\left( a_2 \leq \varepsilon\ | X_2, \dots, X_n \right) 
&\leq& \widetilde{C} \left( \frac{1}{r\sqrt{1 - u}} \left[\varepsilon + \frac{1}{\text{LCD}_{c\sqrt{n}, r}(Y)}\right] +  e^{-2c^2(1 - u)n}\right)\nonumber \\
&\leq& C \left(\varepsilon + \frac{1}{\sqrt{n}} + e^{-c_1n} \right),\label{a2 estimate}
\end{eqnarray}
where $c_1, C, \widetilde{C} > 0$ are absolute constants. 
Note that in the above expression all $\left(X_j\right)_{j = 2}^n$ are fixed and the probability is taken with respect to the random vector $X_1.$

Now we unfix all random vectors $X_2, \dots, X_n.$  Then,
\begin{eqnarray*}
	\mathbb{P}\left(a_2 \leq \varepsilon \text{ or } b_2 \geq t  \right) &=& \mathbb{E}_{X_2, \dots, X_n} \mathbb{P}_{X_1}\left( a_2 \leq \varepsilon \text{ or } b_2 \geq t\right) \\ 
	&=& \mathbb{E}_{X_2, \dots, X_n} \rchi_{\mathcal{E}}\mathbb{P}_{X_1}(a_2 \leq \varepsilon \text{ or } b_2 \geq t  ) + \mathbb{E}_{X_2, \dots, X_n} \rchi_{\mathcal{E}^c}\mathbb{P}_{X_1}(a_2 \leq \varepsilon \text{ or } b_2 \geq t  ) \\
	&\leq&  \mathbb{P}(a_2 \leq \varepsilon\ | X_2, \dots, X_n) + \mathbb{P}(\mathcal{E}^c).
\end{eqnarray*}
Combining the probability estimates in~\eqref{event} and~\eqref{a2 estimate}, we get
\begin{equation*}
\mathbb{P}\left(a_2 \leq \varepsilon \text{ or } b_2 \geq t  \right) \leq C\left(\varepsilon + \frac1{\sqrt{n}} + e^{-c_1n}\right) + \left(\frac2{t^2} + 5u^{n - 2}\right).
\end{equation*}
Repeating this argument for $a_k$ and $b_k$ for $k = 3, \dots, n, $ we obtain for any $\varepsilon,\,t > 0$ and $u \in (0, 1)$
\begin{eqnarray}
\mathbb{P}\left(\frac{a_k}{b_k} \leq \frac{\varepsilon}{t} \right) &\leq& C\left(\varepsilon + \frac1{\sqrt{n}} + e^{-c_1n} + \frac2{t^2} + 5u^{n - 2}\right)  \leq C_1\left(\varepsilon + \frac1{\sqrt{n}} + \frac1{t^2} \right) \label{ak and bk bound}
\end{eqnarray}
where $C, C_1, c_1 > 0$ are absolute constants.

Now we proceed to estimate the sum of $\left(\frac{a_k}{b_k}\right)^2$ in~\eqref{sum}:  
\begin{eqnarray*}
	\mathbb{P}\left(\norm{A^{-1}x} \leq \frac{\varepsilon}{t} \sqrt{n} \right)  &\leq&  
	\mathbb{P}\left(\frac1n \sum_{k = 2}^{n} \left( \frac{a_k}{b_k} \right)^2 \leq \frac{\varepsilon^2}{t^2} \right)  \nonumber \\ 
	&\leq& \mathbb{P}\Big(\exists\, k_1, \dots, k_{\left\lfloor\frac n2\right\rfloor} \in \{2, \dots, n\} \text{ such that } \\ &&\left( \frac{a_{k_i}}{b_{k_i}} \right)^2 \leq 2 \frac{\varepsilon^2}{t^2}  
	\text{ for all } i\leq \left\lfloor\frac n2\right\rfloor \Big)\\ 
	&=& \mathbb{P}\left(\sum_{k=2}^{n} \rchi_k \geq \left\lfloor\frac n2\right\rfloor   \right)
	\ \leq \ \frac2n \sum_{k = 2}^{n} \mathbb{P}\left(  \left( \frac{a_k}{b_k} \right)^2 \leq 2\frac{\varepsilon^2}{t^2}   \right), \nonumber
\end{eqnarray*}
where we denoted by $\rchi_k$ the indicator function of the event $\mathcal{E}_k = \left\{ \left( \frac{a_k}{b_k} \right)^2 \leq 2\frac{\varepsilon^2}{t^2} \right\}$ and in the last step  used Markov's inequality. Using the bound in~\eqref{ak and bk bound}, we finally obtain 
\begin{equation*}
\mathbb{P}\left(\norm{A^{-1}x} \leq \frac{\varepsilon}{t} \sqrt{n} \right) \leq 2 C_1\left(2\varepsilon + \frac1{\sqrt{n}} + \frac1{t^2}\right).
\end{equation*}
Together with an estimate in~\eqref{expectation of x}, we have
\begin{eqnarray*}
	\mathbb{P}\left( s_n(A) \leq \frac{\tau t}{\varepsilon}  n^{-\frac12}   \right) &\geq& \mathbb{P}\left( \norm{x} \leq \tau,\   \norm{A^{-1}x} \geq \frac{\varepsilon}{t} \sqrt{n}   \right) \\
	&\geq& 1 -  C_2\left(\varepsilon + \frac1{\sqrt{n}} + \frac1{t^2}  + \frac1{\tau^2} \right).
\end{eqnarray*}
Since the above statement holds for arbitrary $\varepsilon, t, \tau >0,$ the choice $t = \tau = \frac1{\sqrt{\varepsilon}}$ gives the desired quantitative estimate in Theorem~\ref{main theorem}.
\hfill$\scriptstyle\blacksquare$

\bigskip

\noindent{\bf Aknowledgements.} The author thanks Nicole Tomczak-Jaegermann and Alexander Litvak for their constant support and very valuable suggestions and references. The author is also grateful to Konstantin Tikhomirov for inspiring conversations and encouragement.
A part of this work was done when the author participated in the program ``Geometric Functional Analysis and Applications" in the Fall 2017 at the Mathematical Science Research Institute in Berkeley, California, USA, and was partially supported by the National Science Foundation under Grant No.~DMS-$1440140.$

\address

\end{document}